\documentclass[11pt,reqno]{amsart}

\usepackage{amssymb,amsmath}

\newcommand{\diam}{\mbox{\rm diam}}
\newcommand{\be}{\begin{eqnarray}}
\newcommand{\ee}{\end{eqnarray}}

\newcommand{\card}{\#}
\newcommand{\const}{\mbox{\rm const}}

\newcommand{\es}{\emptyset}
\newcommand{\half}{\frac{1}{2}}

\newcommand{\eps}{{\mbox{$\epsilon$}}}

\newcommand{\Th}{{\theta}}

\newcommand{\R}{{\mathbb R}}

\newcommand{\Nat}{{\mathbb N}}

\newcommand{\Rk}{{\mathcal R}}

\newcommand{\Hau}{{\mathcal H}}

\newcommand{\E}{{\bf E\,}}

\newcommand{\Pro}{{\bf P}}
\newcommand{\Fk}{{\mathcal F}}
\newcommand{\Gk}{{\mathcal G}}

\newcommand{\bn}{{\bf n}}

\newcommand{\Uk}{{\mathcal U}}

\newcommand{\lam}{\lambda}
\newcommand{\om}{\omega}
\newcommand{\gam}{\gamma}

\newcommand{\Leb}{{\mathcal L}}

\def\Exp{{\bf E}}
\def\Var{{\bf Var}}
\def\Dk{{\mathcal D}}
\def\I{\mathcal I}
\newtheorem{theorem}{Theorem}[section]
\newtheorem{lemma}[theorem]{Lemma}
\newtheorem{cor}[theorem]{Corollary}
\newtheorem{prop}[theorem]{Proposition}

\theoremstyle{definition}

\theoremstyle{remark}

\numberwithin{equation}{section}

\input epsf.sty

\begin{document}

\thispagestyle{empty}

\title[projections of zero length]
{{ The sharp Hausdorff measure condition for length of projections}}

\author{Yuval Peres}
\address{Yuval Peres,
Department of Statistics, University of California, Berkeley. \newline
 {\tt peres@stat.berkeley.edu}}

\author{Boris Solomyak}\thanks{Research of Peres was
partially supported by NSF grants \#DMS-0104073 and \#DMS-0244479.
Part of this work was done while he was visiting Microsoft Research.
Research of Solomyak was supported in part by NSF grant
 \#DMS-0099814.
}

\address{Boris Solomyak, Box 354350, Department of Mathematics,
University of Washington, Seattle, WA 98195.
{\tt solomyak@math.washington.edu}}

\subjclass{Primary: 28A80. 
Secondary: 28A75, 
           60D05, 
           28A78 
}

\begin{abstract}
In a recent paper, Pertti Mattila asked  which gauge functions $\varphi$ have the
property that
for any Borel set $A \subset \R^2$ with Hausdorff measure $\Hau^\varphi(A) > 0$, the
projection of $A$ to almost every line has positive length.
We show that finiteness of $\int_0^1\frac{\varphi(r)}{r^2} dr$, 
 which is known to be sufficient for this property, is also necessary
for regularly varying $\varphi$.
Our proof is based on a random construction adapted to the gauge function.
\end{abstract}

\maketitle

\section{{\bf Introduction}} \label{sec:intro}

A classical theorem of Marstrand \cite{mars} 
states that if a planar Borel set $A$ has Hausdorff dimension
strictly greater than $1$, then the orthogonal projection of $A$
to almost every line has positive length.
If $A$ has dimension 1, the situation is more delicate.
Recall that given a positive function $\varphi$ on $(0,\infty)$,
 the  Hausdorff measure
$\Hau^\varphi$ is defined
$
\Hau^\varphi(A)= \lim_{\eps  \downarrow 0} \inf \Bigl\{ \sum_i \varphi(\diam \: A_i) \, : \,
A \subset \bigcup_i A_i \, ,  \: \diam \:  A_i < \epsilon \Bigr\}.
$

In his definitive survey on Hausdorff dimension and projections,
Mattila \cite{mattila} asked  which gauge functions $\varphi$ have the property that
for any Borel set $A \subset \R^2$ with $\Hau^\varphi(A) > 0$, the projection
of $A$ to almost every line has positive length.
In this paper we settle Mattila's question, showing that an integral
 condition known to be sufficient for this property is also necessary; a partial result
 in this direction was  obtained in \cite{JM}. The solution also clarifies the relation
 between Hausdorff measures and  integral-geometric measure.

\medskip

\noindent{\bf Notation.} Let $\eta$ be the
isometry-invariant measure on the space of all lines in $\R^2$, and
 define the integral-geometric measure $\I^1$ on Borel sets in $\R^2$
by $\I^1(A)= \int \#(A \cap \ell) \, d\eta(\ell)$.
(See (\ref{integg}) and \cite{matbook}, Section 5.14.)
Let $p_\theta$ denote  the
 orthogonal projection from $\R^2$ onto the line through the origin making angle
 $\theta$ with the horizontal axis. 
We denote by $\Leb^m$ the $m$-dimensional Lebesgue measure.
We write $f \asymp g$ if $f \le \const\cdot g$ and $g \le \const \cdot f$ 
for some
uniform constant.

\begin{theorem} \label{th-main} Let $\varphi:(0,\infty) \to  (0,\infty)$ be a
(weakly) increasing
function such that
\begin{equation} \label{regul}
\varphi(r)/r^2 \mbox{ \it is (weakly) decreasing.}
\end{equation}
Then the following are equivalent:
\item{\bf (i)} $\int_0^1 \frac{\varphi(r)}{r^2}\,dr < \infty$.
\item{\bf (ii)} If a Borel set $A\subset \R^2$ satisfies $\Hau^\varphi(A) > 0$,
then 
\begin{equation} \label{faveq}
\Leb^1(p_\theta(A))>0 \mbox{ \rm for almost all  } \theta \in [0,\pi) \, .
\end{equation}
\item{\bf (iii)} $\Hau^\varphi$ is absolutely continuous to  $\I^1$ on
Borel sets in $\R^2$.
\end{theorem}
The implication (i)$\Rightarrow$(ii) is known:
Indeed, (i) gives $\int_0^1 \frac{1}{r}\, d\varphi(r) < \infty$ using integration by parts.
 By \cite{Ca}, Theorem IV.1,  this implies that for sets $A$ with
$\Hau^\varphi(A) > 0$, the one-dimensional capacity of $A$ is
positive; a theorem of Kaufman \cite{kauf}  (see also \cite{matbook}, Cor.\ 9.8)
then yields (\ref{faveq}). Clearly (ii)$\Rightarrow$(iii) by definition.
We establish the new implication (iii)$\Rightarrow$(i) by a
random construction; it is intriguing that so far, deterministic constructions
yield less sharp results. See Corollary \ref{cor-int} below
for an extension of Theorem \ref{th-main} to higher dimensions.
Note that for a Borel set $A \subset \R^2$, having positive
one-dimensional capacity is not necessary for (\ref{faveq});
rectifiable curves have zero capacity but certainly satisfy
(\ref{faveq}).

What if the gauge function does not satisfy the regularity condition (\ref{regul})?
In the following lemma we show that if $\varphi$ is increasing, then
we can find another gauge function $\varphi_1$, so that
$\varphi_1(r)/r^2 \downarrow$ and $\Hau^{\varphi_1} \asymp \Hau^\varphi$. (Note, however,
that it could happen that $\int_0^1 \frac{\varphi(r)}{r^2}\,dr=+\infty$, but
$\int_0^1 \frac{\varphi_1(r)}{r^2}\,dr < \infty$; in that case,
for $A$ with $\Hau^\varphi(A)>0$, almost every projection
 has positive length.)

\begin{lemma} \label{lem-gauge} Suppose that $\varphi:(0,\infty) \to (0,\infty)$
is weakly increasing. Then there exists a weakly increasing
$\varphi_1:(0,\infty) \to (0,\infty)$
such that  $\varphi_1(r)/r^2 \downarrow$, and
\begin{equation}\label{eq-gauge}
\Hau^{\varphi_1}(A) \le \Hau^\varphi(A) \le C\cdot \Hau^{\varphi_1}(A)
\end{equation}
for all Borel sets $A \subset \R^2$.
\end{lemma}

\noindent{\bf Remark.} Joyce and M\"orters \cite{JM} constructed a
set $A \subset \R^2$ such that
$0 < \Hau^\phi(A) < \infty$, for some $\phi$ satisfying
$\int_0^1 \frac{\phi(r)^\alpha}{r^{1+\alpha}}\,dr < \infty$ for every $\alpha>1$,
such that $\Leb^1(p_\theta(A))=0$ for {\em all} $\theta$.

Next, we discuss higher dimensions.
Let $\Hau^m$ denote $m$-dimensional Hausdorff measure in $\R^n$.
Recall from \cite{matbook}, Section 5.14,
that the integral-geometric measure $\I^m$ is defined on Borel sets
$A \subset \R^n$ by
\begin{equation} \label{integg}
\I^m(A)= \int \! \int \#\Bigl(A \cap p_V^{-1}(a)\Bigr) \, d\Hau^m(a) \, d\gamma_{n,m}(V) \,
\end{equation}
where $p_V$ is the orthogonal projection to $V$, and $\gamma_{n,m}$ is the
isometry-invariant measure on the Grassman manifold of $m$-dimensional subspaces
of $\R^n$. The next corollary generalizes Theorem \ref{th-main}.
\begin{cor} \label{cor-int}
 Let $m<n$  and suppose  $\psi$ is a  positive
function on $(0,\infty)$ such that $\psi(r)/r^{m-1}$ is weakly  increasing,
$\psi(r)/r^{m+1}$ is weakly decreasing (these regularity conditions could be relaxed).
Then the following are equivalent:
\item{\bf (i)} $\int_0^1 \frac{\psi(r)}{r^{m+1}}\,dr < \infty$.
\item{\bf (ii)} If a Borel set $\Lambda \subset \R^n$ satisfies $\Hau^\psi(A) > 0$,
then $\Leb^m(p_V(\Lambda))>0$ for $\gamma_{n,m}$ almost all $m$-dimensional subspaces
$V \subset \R^n$.
\item{\bf (iii)} $\Hau^\psi$ is absolutely continuous to  $\I^m$ on
Borel sets in $\R^n$.
\end{cor}
The implication (i)$\Rightarrow$(ii) is known:
it follows from combining \cite{Ca}, Theorem IV.1,  with \cite{matbook}, Cor.\ 9.8.
 Since (ii)$\Rightarrow$(iii) is obvious, we will only need to prove
 (iii)$\Rightarrow$(i), and this will follow from the same construction we
 use to establish the corresponding implication in Theorem \ref{th-main}.

\section{{\bf Outline of the construction}}
To prove the implication (iii)$\Rightarrow$(i) in Theorem \ref{th-main},
we establish the following.

\begin{prop} \label{prop-main}
 Let $\varphi$ be a weakly increasing function that satisfies (\ref{regul})
and $\int_0^1 \frac{\varphi(r)}{r^2} dr = \infty$. Then there
exists a  Borel set $A\subset \R^2$ such that $\Hau^\varphi(A) > 0$
and $\Leb^1(p_\theta(A))=0$ for almost all $\theta \in [0,\pi)$.
\end{prop}
The set $A$ is constructed as a random Cantor set in the plane; we will show
that it has all the desired properties almost surely.
This set is a (substantial)
modification of a stochastically self-similar set constructed
in \cite{buffon}.

Denote by $\Gk_k$ the collection of $4^k$ (closed) dyadic
subsquares of the unit square $[0,1]^2$ having  side length
$2^{-k}$. We consider all dyadic subsquares as a rooted tree,
with $[0,1]^2$ being the root and $\Gk_k$ being the set of nodes
at the $k$th level. For each node there are four edges leading to
nodes at the next level (its ``offspring'').

First we inductively define random subsets $\Fk_k \subset \Gk_k$ for $k\ge 0$,
and let $\Rk_k := \bigcup \{Q:\ Q\in \Fk_k\}$.
We start with $\Fk_0 = \Gk_0 = \{[0,1]^2\}$.
Passing from $\Fk_k$ to $\Fk_{k+1}$ is either deterministic or random,
depending on $k$.
Given $\Fk_{k}$, a set  of
dyadic squares of side length $2^{-k}$, we partition each of them into
four dyadic subsquares of side length $2^{-k-1}$.
Now, if the step is deterministic, then we keep all the squares, so that
$\#\Fk_{k+1} =4 \# \Fk_k$ and $\Rk_{k+1} = \Rk_k$.
If the step is random, then we choose one of the four
subsquares uniformly at random, all these choices being independent.
This way we obtain a subset $\Fk_{k+1}$ with the same cardinality as
$\Fk_k$.  The gauge function $\varphi$ dictates whether we make a deterministic
or random step. If we
do $\alpha_k$ deterministic steps for $i< k$, then $\card \Fk_k =
4^{\alpha_k}$, and  we make sure that $4^{-\alpha_k} \asymp \varphi(2^{-k})$.
Once $\Fk_k$ and $\Rk_k$ are defined, we consider $\Rk = \bigcap_{k=1}^\infty
\Rk_k$, which satisfies $0 < \Hau^\varphi(\Rk) < \infty$.

We do not know whether the resulting random set has a.e.\
projection of zero length (almost surely),
so we modify the construction slightly,
removing certain squares at specific levels chosen in advance. This
yields a new random set $\Rk'$ which has all the desired properties.


\section{{\bf Construction in detail}}

Consider a gauge function $\varphi(r)$ as in Theorem~\ref{th-main}.
Without loss of generality, we can assume that $\varphi(1)=1$.
Let
$$
\alpha_n = \lfloor \log(\varphi(2^{-n}))/\log(1/4) \rfloor\ \ \ \mbox{for}\ \
n\ge 0,
$$
so that $\alpha_0=0$ and
\begin{equation} \label{defal}
\varphi(2^{-n}) \le 4^{-\alpha_n} \le 4\cdot \varphi(2^{-n}).
\end{equation}
Since $\varphi(r)$ is increasing, we have $\alpha_n \ge \alpha_{n-1}$, and since
$\varphi(r)/r^2$ is decreasing, we have $\alpha_n \le \alpha_{n-1}+1$.
Now,
$$
\infty = \int_0^1\frac{\varphi(r)}{r^2}dr = \sum_{n=0}^\infty
\int_{2^{-n-1}}^{2^{-n}} \frac{\varphi(r)}{r^2}dr  \le \sum_{n=0}^\infty
2^{n+1} \varphi(2^{-n}) \le \sum_{n=0}^\infty 2^{n+1} 4^{-\alpha_n},
$$
so
\begin{equation} \label{series}
\sum_{n=0}^\infty \lam_n = \infty,\ \ \ \mbox{where}\ \ \lam_n:=
2^n 4^{-\alpha_n}.
\end{equation}

Now we specify when we perform deterministic steps. The transition from
$\Fk_{n-1}$ to $\Fk_n$ is deterministic if $\alpha_n>\alpha_{n-1}$,
that is, $\alpha_n = \alpha_{n-1}+1$. Otherwise, the
transition is random, as described in Section 2. This is well-defined,
and we make $\alpha_n$ deterministic steps on levels $i<n$, so
$\# \Fk_n = 4^{\alpha_n}$. Recall that $\Rk_n = \bigcup_{B\in \Fk_n} B$
and $\Rk = \bigcap_{n=1}^\infty \Rk_n$.
Consider the probability measure $\mu$ on $\Rk$ defined as the weak$^*$
limit of uniform measures on $\Rk_n$. We claim that there exists $C_1>0$
such that
\begin{equation} \label{est1}
\mu(B_r) \le C_1\cdot \varphi(r)
\end{equation}
for any ball of radius $r$. Indeed,
up to a multiplicative constant, we can replace the ball $B_r$ in
(\ref{est1}) by a square $B \in \Fk_{n}$ where $2^{-n} \asymp r$.
We have $\mu(B) = 1/\#\Fk_n = 4^{-\alpha_n} \asymp \varphi(2^{-n})$,
and the desired estimate follows. Thus, $\Hau^\varphi(\Rk)>0$ by the Mass
Distribution Principle. We can assume that $\lam_n$ are bounded,
\begin{equation} \label{supcon}
\lam_n < C_2,
\end{equation}
for all $n$.
Indeed, the set $\Rk_n$ consists of $4^{\alpha_n}$ squares of side length
$2^{-n}$. If (\ref{supcon}) does not hold, then  $\inf (4^{\alpha_n} 2^{-n})
=0$, hence $\Hau^1(\Rk)=0$, and then all projections of $\Rk$ have
zero length.

\medskip

Next we specify how we modify the set $\Rk$ by removing certain squares
at  prescribed levels.
We label the four dyadic subsquares of a square as in Figure 1.
\begin{figure}
\centerline{\epsfxsize=1.in \epsffile{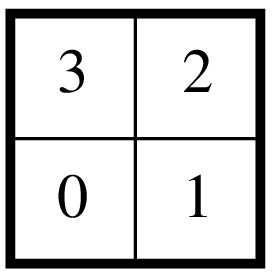}}
\caption{Labeling subsquares.}
\end{figure}
This labeling induces a natural addressing scheme for each dyadic
square $B\in \Gk_k$. The address has length $k$ and the symbols
are from $\{0,1,2,3\}$; we write it as $\om(B) = \{\om_i(B)\}_{i=1}^{k}$.

For $B \in \Gk_{n}$ and $1 \le i \le n$
define
$$
Y_i(B) =
\left\{ \begin{array}{ll} \lam_i,
                                      & \mbox{if}\ \om_{i}(B)\in \{0,2\}; \\
                                   0, & \mbox{if}\ \om_{i}(B)\in \{1,3\}.
\end{array} \right.
$$
We can consider $\Gk_{n}$ as a discrete probability space with the
uniform measure. Then $\{Y_i\}_{i=1}^{n}$ are independent random variables,
such that $\Pro(Y_i = \lam_i) = \Pro(Y_i = 0) = \half$.
Let
$$
S_n(B) = \sum_{i=1}^{n} Y_{k_i}(B),
$$
where $k_i$ is the $i$-th deterministic step, so that $\alpha_{k_i-1}=i-1$ and
$\alpha_{k_i}=i$. Denote
$$
\gam_i := \lam_{k_i} = 2^{k_i} 4^{-\alpha_{k_i}} = 2^{k_i} 4^{-i};
$$
then
$\Exp (S_n) = \half \sum_{i=1}^{n} \gam_{i}$ and
$\Var (S_n) =  \frac{1}{4} \sum_{i=1}^{n} \gam_i^2$.
Observe that $S_n(B)$ depends only on the digits of $B$ corresponding to
deterministic steps, hence it is independent of the event  $B\in \Fk_{k_n}$
which can be considered on the same probability space $\Gk_{k_n}$.

We say that $B\in \Gk_{k_n}$ is {\it deviant\/} if
$$
|S_n(B) - \Exp(S_n)| >  \Exp(S_n)/2.
$$
Using Chebyshev's inequality we obtain
$$
\Pro(B \mbox{\ is  deviant}) \le \frac{\Var(S_n)}{(\half \Exp(S_n))^2}
= \frac{\sum_{i=1}^{n} \gam_{i}^2}{(\sum_{i=1}^{n} \gam_{i})^2} \le
C_2\left(\sum_{i=1}^{n} \gam_{i}\right)^{-1},
$$
since $\gam_{i}^2 \le C_2 \gam_{i}$ by (\ref{supcon}).
Denoting by $\Dk_{n}$ the collection of deviant squares in $\Gk_{k_n}$
and using independence of $B$ being deviant from the event
$(B\in \Fk_{k_n})$,
we obtain that
$$
\# (\Fk_{k_n} \cap \Dk_{n})
\le C_2 \cdot \# \Fk_{k_n} \left(\sum_{i=1}^{n} \gam_{i}\right)^{-1}\,.
$$
By the definition of measure $\mu$ we have
\begin{equation} \label{dop2}
\mu\left(\Rk \setminus \bigcup_{B \in \Dk_{n}} B \right) \ge 1 -
C_2 \left(\sum_{i=1}^{n} \gam_{i}\right)^{-1}\,.
\end{equation}

Note that $\alpha_i = j$ for $k_j \le i < k_{j+1}$, hence denoting $k_0=0$
we have
$$
\sum_{i=0}^{k_n-1} \lam_i =
\sum_{j=0}^{n-1} \sum_{i=k_j}^{k_{j+1}-1} \lam_i
= \sum_{j=0}^{n-1} 4^{-j} \sum_{i=k_j}^{k_{j+1}-1} 2^i <
\sum_{j=0}^{n-1} 4^{-j} 2^{k_{j+1}} = 4 \sum_{j=1}^n \gam_j.
$$
Observe that $k_n\to \infty$, otherwise $\Rk$ is a finite set.
Thus, $\sum_{i=1}^{n} \gam_{i}\to \infty$ as $n\to\infty$ by
(\ref{series}), and
we can find $n(j) \in \Nat$, $j\ge 1$,  such that
\begin{equation} \label{dop3}
C_2^{-1}\sum_{i=1}^{n(j)} \gam_{i} > 2^{j+1}.
\end{equation}
Let
$$
\Rk' := \Rk \setminus \bigcup_{j=1}^\infty \bigcup_{B\in \Dk_{n(j)}} B.
$$
We have
$\mu(\Rk') \ge 1 - \half \sum_{j=1}^\infty 2^{-j} = \half$ by (\ref{dop2})
and (\ref{dop3}).  Thus, $\mu|_{\Rk'}$ is a positive measure, and
$\mu(B_r \cap \Rk') \le \mu(B_r) \le C_1\cdot \varphi(r)$ for any ball $B_r$ of
radius $r$, hence $\Hau^\varphi(\Rk') > 0$ by the Mass Distribution Principle.


\section{{\bf Proof of Proposition \ref{prop-main} and  Theorem~\ref{th-main} }}

Denote
$$
\Rk_{k_n}' = \Rk_{k_n} \setminus \bigcup_{B\in \Dk_n} B.
$$

\begin{lemma} \label{lem-main}
Let $\ell$ be a line intersecting $[0,1]^2$
that does not hit any vertices of the squares in $\Gk_{k_n}$. Then
$$
\Pro(\Rk'_{k_n} \cap \ell \ne \es) \le 64\cdot
\left(\sum_{i=1}^{n} \gam_i\right)^{-1}.
$$
\end{lemma}

\noindent {\sc Proof of Theorem~\ref{th-main} assuming
Lemma~\ref{lem-main}.}
Observe that $\Rk' = \bigcap_{j=1}^\infty \Rk_{k_{n(j)}}'$, hence
Lemma~\ref{lem-main} implies that the probability of $\ell$ hitting
$\Rk'$ equals zero.
Let $\Th \in [0,\pi]$ be such that the
line $y\cos\Th = x\sin\Th$ is orthogonal to $\ell$, and let $\bn$
be the unit normal vector for $\ell$.
Then by Fubini's Theorem,
$$
\Exp\left[ \Leb^1(p_\theta(\Rk')) \right] = \int_{\R} \Pro(\Rk' \cap
(\ell + t \bn) \ne \es)\,dt =0 \,.
$$
The proposition, and hence  the theorem, follow  by integrating over $\Th$. \qed

\medskip

\noindent {\sc Proof of Lemma~\ref{lem-main}.}
Let $\alpha$ be the angle that $\ell$ forms with the horizontal.
First we assume that $\alpha \in [0,\pi/2]$ and
then indicate how to consider the case $\alpha \in (\pi/2,\pi]$.
Let
$$
A_n(\ell) = \card\{B\in \Gk_{k_n}:\ B \cap \ell \ne \emptyset\}.
$$
Observe that
\begin{equation} \label{eq-estan} A_{n}(\ell) \le  2^{k_n+1}.
\end{equation}
To verify this  we may assume, using symmetry, that $\ell$ forms
an angle in $[0,\pi/4]$ with the horizontal. Then $\ell$
intersects at most two squares in each of the $2^{k_n}$ columns of
$\Gk_{k_n}$, and (\ref{eq-estan}) follows.

Recall that none of the squares in $\Rk_{k_n}'$ are
deviant, i.e.
$$
B \in \Rk_{k_n}' \ \Rightarrow\ (1/4) \sum_{i=1}^{n} \gam_i \le S_n(B) \le
(3/4) \sum_{i=1}^{n} \gam_i,
$$
where $S_n(B) = \sum \{\gam_i:\ i\le n,\,\om_{k_i}(B)\in \{0,2\}\}$.

(i) Say that $B \in \Gk_{2n}$ is {\bf $0$-rich} if
$$
\sum\{\gam_i:\ i \le n,\,\om_{k_i}(B) = 0\} \ge (1/8)
\sum_{i=1}^{n} \gam_i.
$$

(ii) Say that $B \in \Gk_{2n}$ is {\bf $2$-rich} if
$$
\sum\{\gam_i:\ i \le n,\,\om_{k_i}(B) = 2\} \ge (1/8)
\sum_{i=1}^{n} \gam_i.
$$

\noindent By the definition of deviant
squares, every non-deviant square is
either  $0$-rich or  $2$-rich  (or both).
Consider the events
$$
Z_i = \Big\{ \exists\, B\in \Fk_{k_n}:\ B\ \mbox{is $i$-rich}\ \&\ B\cap
\ell \ne \emptyset \Big\} \ \ \ \mbox{for}\ \ i=0,2.
$$
Since $\Fk_{k_n}'$ contains only  non-deviant
squares from $\Fk_{k_n}$, we have
$$
\Pro(\Rk_{k_n}' \cap \ell \ne \es) \le \Pro(Z_0)+\Pro(Z_2).
$$
Let us estimate $\Pro(Z_0)$. We have
\begin{equation} \label{eq-condex}
\E[\card\{B\in \Fk_{k_n}:\ B\cap \ell \ne \es\}\,|\, Z_0] \le
\frac{\E[\card\{B\in \Fk_{k_n}:\ B\cap \ell \ne \es\}]}{\Pro(Z_0)}\,.
\end{equation}
Observe that
\begin{equation} \label{eq-dop1}
\Pro(B\in \Fk_{k_n}) = \#\Fk_{k_n}/\#\Gk_{k_n} = 4^{n}/2^{2k_n}=
4^{n-k_n}
\end{equation}
for any $B\in \Gk_{k_n}$.
Writing
$$
\card\{B\in \Fk_{k_n}:\ B\cap \ell \ne \es\} =
\sum_{B\in \Gk_{k_n}} {\bf 1}_{\{B\in \Fk_{k_n}:\ B\cap \ell \ne \es\}}
$$
and using (\ref{eq-dop1}) we obtain by (\ref{eq-estan}) that
\begin{equation} \label{eq-ex}
\E\Big[\card\{B\in \Fk_{k_n}:\ B\cap \ell \ne \es\}\Big] =
A_n(\ell)\cdot 4^{n-k_n} \le 2\cdot 2^{2n-k_n}.
\end{equation}
Thus by (\ref{eq-condex}),
\begin{equation} \label{eq-condez}
{\Pro(Z_0)} \le
\frac{2^{2n-k_n+1}} {\E[\card\{B\in \Fk_{k_n}:\ B\cap \ell \ne \es\}\,|\, Z_0]}
\, .
\end{equation}
 It remains to estimate the denominator in (\ref{eq-condez}) from below. Let
$$
\Psi_0 := \{Q\in \Gk_{k_n}:\ Q\ \mbox{is $0$-rich}\ \&\ Q\cap \ell \ne \es\}.
$$
Order the squares in $\Gk_{k_n}$ hit by $\ell$ from left to right and from
bottom to top. This is a total order by the assumption on slope of the
line $\ell$. For $Q\in \Psi_0$ consider the event
$$
Y_Q = \Big\{ \mbox{ $Q$\ \  is the first square in $\Psi_0$ such that
 $Q\in \Fk_{k_n}$} \Big\}.
$$
Then $Z_0 = \bigcup_{Q\in \Psi_0} Y_Q$ is a disjoint union, and
so, for any random variable $f$,
\begin{equation} \label{eq-conv}
\E[f\,|\,Z_0] = \sum_{Q\in \Psi_0} \frac{\Pro(Y_Q)}{\Pro(Z_0)}
\E[f\,|\,Y_Q] \ge \min_{Q \in \Psi_0}\E[f\,|\,Y_Q].
\end{equation}
Fix $Q\in \Psi_0$. We have
\begin{equation} \label{eq-condex2}
\E[\card\{B\in \Fk_{k_n}:\ B\cap\ell\ne\es\}\,|\,Y_Q] =
\sum_{B\in \Gk_{k_n}:\ B\cap\ell\ne\es} \Pro(B\in \Fk_{k_n}\,|\,Y_Q).
\end{equation}
Fix $i$ such that $\om_{k_i}(Q)=0$, and denote by
$\widetilde{Q}$ the dyadic square in $\Gk_{k_i-1}$ having the
address $\om(\widetilde{Q}) = \om_1(Q)\ldots w_{k_i-1}(Q)$. The fact
that $Q\in \Fk_{k_n}$ implies that
$\widetilde{Q}\in \Fk_{k_i-1}$.
(Note that by definition, $[0,1]^2
\supset \Rk_{k_1} \supset \ldots \supset \Rk_{k_n}$.)
Recall that the step $k_i$ is deterministic, so all four offspring of
$\widetilde{Q}$ are in $\Fk_{k_i}$.
Since the slope of $\ell$ is positive, $\ell$
intersects at least $\half 2^{k_n-k_i}$ squares $B\in \Gk_{k_n}$ whose
addresses start with $\om(\widetilde{Q})k$, for $k\in \{1,2,3\}$
(see Figure 1). For each of these squares we have (using the
independence of $Y_Q$ from the random choices involving the
descendants of $\om(\widetilde{Q})k$ with  $k\in \{1,2,3\}$), that
$$
\Pro(B\in \Fk_{k_n}\,|\,Y_Q) = \Pro(B\in \Fk_{k_n}\,|\,
\widetilde{Q}\in \Fk_{k_i-1}) = 4^{k_i-k_n+n-i}.
$$
This is because we made $(n-i)$ deterministic steps between
stages $k_i$ and $k_n$, and hence $(k_n-k_i)-(n-i)$ random steps.
Therefore, the sum of $\Pro(B\in \Fk_{k_n}\,|\,Y_Q)$ over the set
of squares
$$
\Uk_i=\Big\{ B\in \Gk_{k_n}, \, : \ B\cap\ell\ne\es,\ \
\{\om_j(B)\}_1^{k_i-1} = \{\om_j(Q)\}_1^{k_i-1},\ \ \om_{k_i}(B) \in
\{1,2,3\} \Big\},
$$
is at least $\half 2^{k_n-k_i}\cdot 4^{k_i-k_n+n-i}= \half 2^{2n-k_n}\gam_i$.
Notice that the sets $\Uk_i$
are disjoint for distinct $i$ with
$\om_{k_i}(Q)=0$. Thus, by the definition of $0$-rich squares,
the right-hand side of (\ref{eq-condex2}) is at least
$$
(1/2) 2^{2n-k_n} \cdot \sum \{\gam_i:\ i\le n,\ \om_{k_i}(B)=0 \} \ge (1/16)
2^{2n-k_n} \cdot \sum_{i=1}^{n} \gam_i \,.
$$
Therefore by (\ref{eq-conv}),
$$
\E[\card\{B\in \Fk_{k_n}:\ B\cap \ell \ne \es\}\,|\, Z_0] \ge 2^{2n-k_n-4}
\cdot \sum_{i=1}^{n} \gam_i \,.
$$
Together with (\ref{eq-condez}), this
yields that
$$\Pro(Z_0) \le  32 \cdot \left({\sum_{i=1}^{n} \gam_i}\right)^{-1}.
$$

The estimate $P(Z_2) \le 32 \cdot \left({\sum_{i=1}^{n} \gam_i}\right)^{-1}$
is proved similarly. We consider the set $\Psi_2$ of squares in $\Gk_{k_n}$
hit by $\ell$ that are $2$-rich and condition on the
{\em last} square in $\Psi_2$ which belongs to $\Fk_{k_n}$.

This concludes the proof of
the lemma when the angle $\alpha$ is in $[0,\pi/2]$. In the case when
$\alpha \in [\pi/2,\pi]$ we interchange the roles of
the subsquares  $\{0,2\}$ and $\{1,3\}$ and use that for a  non-deviant
square $B$ we have
$$
\sum \{\gam_i:\ i \le n,\ \om_{k_i}(B) \in \{1,3\}\} \ge
(1/4) \sum_{i=1}^{n} \gam_i.
$$
\qed


\section{{\bf Proof of Lemma~\ref{lem-gauge}}}

Let
$$
\varphi_1(s) = \inf_r \left\{(s/r)^2 \varphi(r):\ r\le s \right\}\,.
$$
Then $\varphi_1(s) \le \varphi(s)$, so the left-hand side in (\ref{eq-gauge}) is clear.
For the right-hand side consider a cover of a set $A$ for $\Hau^{\varphi_1}$.
Let $S$ be a set in this cover where $S$ has diameter $s$. Then we can find
$r\le s$ with $(s/r)^2 \varphi(r) < 2 \varphi_1(s)$. Let $v<r$ be the largest number
of the form $s/2^k$. Then $(s/v)^2 \varphi(v) < 8 \varphi_1(s)$. Now $S$ may be
covered by a square of side $s$, hence by $4^{k+1} = 4(s/v)^2$ squares of
side $v/2$. Replacing $S$ by these squares and using them in the Hausdorff
sum for $\Hau^{\varphi}$ we conclude that $\Hau^\varphi(A) \le 32 \Hau^{\varphi_1}(A)$.

The condition $\varphi_1(s)/s^2 \downarrow$ is immediate, since
$\varphi_1(s)/s^2 = \inf \{\varphi(r)/r^2:\ r\le s\}$.
It remains to check that $\varphi_1(s) \uparrow$. Suppose $s< t$. Then
$\varphi_1(t) =\min\{A_1, A_2\}$, where
$$ A_1= \inf_{r\le s}\left\{(t^2/r^2)
\cdot \varphi(r)\right\} \mbox{ \rm and }  A_2=
\inf_{r\in [s,t]}   \left\{(t^2/r^2)\cdot \varphi(r)\right\} \,.
$$
It is clear that $A_1 \ge \varphi_1(s)$, and
$A_2 \ge (t^2/r^2)\cdot
\varphi(s) \ge \varphi(s) \ge \varphi_1(s)$ since $\varphi(r)$ is increasing.
The lemma is proved. \qed

\section{{\bf Proof of Corollary~\ref{cor-int} }}
We only need to prove that (iii)$\Rightarrow$(i). We prove the contrapositive,
i.e., given $\psi$ satisfying the regularity conditions, 
we show that $\Hau^\psi$ is not absolutely continuous to  $\I^m$
as Borel measures on $\R^n$.

Let $\varphi(r)=\psi(r)/r^{m-1}$. Then Proposition \ref{prop-main}
yields a  Borel set $A\subset \R^2$ such that $\Hau^\varphi(A) > 0$
and $\I^1(A)=0$. Define $B \subset \R^n$ as the Cartesian product
$B=A \times [0,1]^{m-1} \times \{(0,\ldots,0)\}$
(appending $n-m-1$ zeros). Then it is standard to verify that
$\Hau^\psi(B)>0$ yet $\I^m(B)=0$. For the latter, observe that for a
generic $m$-dimensional subspace $V$ and $a\in V$, 
the $(n-m)$-dimensional subspace
$p_V^{-1}(a)$ intersected with $\R^{m+1}\times \{(0,\ldots,0)\}$,
is a line which hits $B$ if and only if its projection on the first
two coordinates hits $A$. Thus the inner integral in (\ref{integg}) is
zero for a generic $V$.
\qed

\bibliographystyle{amsplain}

\end{document}